\documentclass[reqno,11pt]{amsart}
\usepackage[letterpaper,margin=1in,footskip=0.25in]{geometry}
\usepackage[utf8]{inputenc}
\usepackage{mathptmx}
\usepackage{microtype}
\usepackage[cal=boondoxupr]{mathalfa}
\usepackage{mathrsfs}
\usepackage{enumerate, setspace, amssymb, amsopn, amsmath, array, pdfsync, url, amsfonts, float,enumitem}
\usepackage{cases}

\usepackage{xcolor}
\definecolor{chianti}{rgb}{0.6,0,0}
\definecolor{meretale}{rgb}{0,0,.6}
\definecolor{leaf}{rgb}{0,.35,0}

\usepackage[colorlinks=true, hyperindex, citecolor=meretale, urlcolor=leaf, linkcolor=chianti]{hyperref}
\usepackage{tikz, tikz-cd}
\usetikzlibrary{plotmarks}
\usepackage{tikz-layers}
\usepackage[capitalise]{cleveref}
\usepackage{booktabs}
\usepackage{multirow}
\usepackage{arydshln}
\definecolor{cerulean}{rgb}{0.0, 0.48, 0.65}
\definecolor{brickred}{rgb}{0.8, 0.25, 0.33}

\usepackage{wrapfig}

\allowdisplaybreaks

\newcounter{rowcntr}[table]
\renewcommand{\therowcntr}{(\arabic{rowcntr})}

\newcolumntype{N}{>{\refstepcounter{rowcntr}\therowcntr}c}

\AtBeginEnvironment{tabular}{\setcounter{rowcntr}{0}}

\makeatletter
\g@addto@macro{\endtabular}{\rowfont{}}
\makeatother
\newcommand{\rowfonttype}{}
\newcommand{\rowfont}[1]{
   \gdef\rowfonttype{#1}#1
}
\newcolumntype{L}{>{\rowfonttype}l}

\newtheorem{Theoremx}{Theorem}
 
\newtheorem{theorem}{Theorem}[section]
\theoremstyle{definition}

\theoremstyle{definition}

\newtheorem{case}{Case}

\theoremstyle{definition}
\newtheorem{definition}[theorem]{Definition}
\newtheorem{remark}[theorem]{Remark}
\theoremstyle{remark}

\newcommand{\Spec}{\operatorname{Spec}}

\usepackage{stmaryrd}

\newcommand{\frk}{\operatorname{frk}}

\newcommand{\N}{\mathbb{N}}

\newcommand{\F}{\mathbb{F}}

\newcommand{\PP}{\mathbb{P}}

\newcommand{\fm}{\mathfrak{m}}
\newcommand{\fp}{\mathfrak{p}}

\newcommand{\fn}{\mathfrak{n}}

\newcommand{\FF}{\mathbb{F}}

\newcommand{\cL}{\mathcal{L}}
\newcommand{\cO}{\mathcal{O}}

\newcommand{\ehk}{\operatorname{e}_{\operatorname{HK}}}
\newcommand{\ehkx}{\operatorname{e}_{\operatorname{HK},X}}

\makeatother

\crefname{Theoremx}{Theorem}{Theorems}
\crefname{maintheorem}{Main Theorem}{Main Theorems}
\crefname{table}{Table}{Tables}
\crefname{figure}{Figure}{Figures}
\crefname{case}{Case}{Cases}

\usepackage[backend=biber, style=alphabetic, datamodel=mrnumber,sorting=anyt,maxalphanames=5,maxbibnames=99]{biblatex}
\usepackage{hyperref}
\usepackage{xurl}
\hypersetup{breaklinks=true}

\DeclareFieldFormat{mrnumber}{
  MR\addcolon\space
  \ifhyperref
    {\href{http://www.ams.org/mathscinet-getitem?mr=#1}{\nolinkurl{#1}}}
    {\nolinkurl{#1}}}

\renewbibmacro*{doi+eprint+url}{
  \iftoggle{bbx:doi}
    {\printfield{doi}}
    {}
  \newunit\newblock
  \printfield{mrnumber}
  \newunit\newblock
  \iftoggle{bbx:eprint}
    {\usebibmacro{eprint}}
    {}
  \newunit\newblock
  \iftoggle{bbx:url}
    {\usebibmacro{url+urldate}}
    {}}

\renewbibmacro{in:}{}

\makeatletter
\renewcommand*{\eqref}[1]{
  \hyperref[{#1}]{\textup{\tagform@{\ref*{#1}}}}
}
\makeatother

\renewcommand{\div}{\mathrm{div}}
\newcommand{\vol}{\mathrm{vol}}

\begin{document}

\title{{H}ilbert--{K}unz multiplicity and \texorpdfstring{$F$}{F}-signature can disagree}

\author[Lee]{Seungsu Lee}
\address{Department of Mathematics, University of Michigan, Ann Arbor, MI 48109 USA}
\email{starlee@umich.edu}

\author[Pande]{Suchitra Pande}
\address{Department of Mathematics, University of Utah, Salt Lake City, UT 84112 USA}
\email{suchitra.pande@utah.edu}

\author[Simpson]{Austyn Simpson}
\thanks{Simpson was supported by NSF postdoctoral fellowship DMS \#2202890.}
\address{Department of Mathematics, Bates College, Lewiston, ME 04240 USA}
\email{asimpson2@bates.edu}

\begin{abstract}
We compute the $F$-signature function of the ample cone of any nontrivial ruled surface over $\mathbb{P}^1_k$ where $k$ is an algebraically closed field of prime characteristic. As an application, we construct a Noetherian $F$-finite strongly $F$-regular ring of prime characteristic admitting two maximal ideals $\mathfrak{n}_1,\mathfrak{n}_2\in \Spec R$ at which the Hilbert--Kunz multiplicity and $F$-signature measure different singularities; that is, $\operatorname{e}_{\operatorname{HK}}(R_{\mathfrak{n}_1})<\operatorname{e}_{\operatorname{HK}}(R_{\mathfrak{n}_2})$ and $s(R_{\mathfrak{n}_1})<s(R_{\mathfrak{n}_2})$. Our calculation of the $F$-signature for the Hirzebruch surfaces also corrects an inaccuracy in a preprint by different authors.
\end{abstract}
\maketitle

\section{Introduction}\label{sec:introduction}

This note concerns two measures of singularity in prime characteristic: the Hilbert--Kunz multiplicity and the $F$-signature. Let $(R,\fm,k)$ be a $d$-dimensional $F$-finite Noetherian local $\FF_p$-algebra. For simplicity in the introduction, we also assume that $R$ is a domain with perfect residue field $k=k^p$. Recall that the \emph{Hilbert--Kunz multiplicity} of $R$, denoted $\ehk(R)$, is defined to be the limit
$$\ehk(R)=\lim\limits_{e\to\infty}\frac{\ell_R(R/\fm^{[p^e]})}{p^{ed}}$$ which was shown to exist in \cite{Mon83}. The \emph{$F$-signature of $R$} on the other hand, denoted $s(R)$, asymptotically counts the splittings of $R\hookrightarrow R^{1/p^e}$ via
$$s(R) = \lim\limits_{e\to \infty}\frac{\frk_R R^{1/p^e}}{p^{ed}},$$ where $$\frk_R R^{1/p^e} = \max\{N\mid \exists R^{1/p^e}\twoheadrightarrow R^{\oplus N}\}$$ is the free rank of the $R$-module $R^{1/p^e}$. This latter limit was shown to exist in \cite{Tuc12}.

Both of these measurements can be thought of as local volumes of a singularity. For example, $\ehk(R)$ and $s(R)$ attain their respective minimal or maximal values of one precisely when $R$ is regular (see \cite[Theorem 1.5]{WY00} and \cite{HL02}). An oft-repeated philosophy concerning these invariants is that values closer to one correspond to less severe singularities, a perspective for which there is ample evidence (see e.g. \cite{AE08,PT18,PS20,CRST21,DS22}). In the present article we raise a counterpoint to this heuristic, inspired by a question posed to the third named author by T. Polstra regarding whether these invariants ever give contradictory verdicts about when one singularity is milder than another.
\begin{Theoremx}\label{theorem:A}
    For every prime integer $p>0$, there exists a Noetherian $F$-finite strongly $F$-regular ring $R$ of characteristic $p$ with maximal ideals $\fn_1,\fn_2\in\Spec(R)$ such that $\ehk(R_{\fn_1})<\ehk(R_{\fn_2})$ and $s(R_{\fn_1})<s(R_{\fn_2})$.
\end{Theoremx}

Simply stated, the inequalities above suggest that $R_{\fn_1}$ is less singular than $R_{\fn_2}$ according to the Hilbert--Kunz multiplicity, whereas the $F$-signature suggests the opposite. The ring $R$ appearing in \cref{theorem:A} is given by the tensor product $R_1\otimes_k R_2$ where $R_1$ and $R_2$ are carefully chosen section rings of the Hirzebruch surface $\mathcal{H}_a$ with parameter $a$ with respect to ample line bundles $\mathcal{L}_1$ and $\mathcal{L}_2$. See \cref{sec:preliminaries} for a summary of the toric geometry of $\mathcal{H}_a$ and for our notational conventions.

The Hilbert--Kunz multiplicity of the section rings of $\mathcal{H}_a$ has been computed already by Trivedi in \cite{Tri16} (see \cref{theorem:Trivedi}). Values for the $F$-signature were claimed in the preprint \cite{HS17}, but unfortunately the calculations there contain an error which prevents their usage in our applications. In essence, the formulas of \emph{op. cit.} do not satisfy the transformation rule for $F$-signature under finite maps (see \cref{rem:transformation-rule} for more details). Therefore, we give a new calculation for the $F$-signature function of the ample cone of $\mathcal{H}_a$ for $a\geq 2$, while the cases of $a=0$ and $a=1$ were handled previously in \cite[Example 4.4.10]{VK12} and \cite[Appendix A]{Lee23} respectively. The majority of the paper is dedicated to this calculation in \cref{sec:calculation}, and our formula (which does satisfy the transformation rule) is summarized below.

\begin{Theoremx}\label{theorem:B}
    Let $a\geq 2$ and let $\pi: X = \mathcal{H}_a \to \PP_k^1$ be a ruled surface over $\PP^1_k$ where $k$ is an algebraically closed field of characteristic $p>0$. Let $D_1$ be the divisor class of a fiber of $\pi$ and $D_4$ be the divisor class of a section of $\pi$. Let $\cL$ be an ample line bundle on $X$ given by $\cL = c D_1 + dD_4$ where $c,d\geq 1$. If $R = \bigoplus_{n\in \mathbb{Z}} H^0(X,\cL^{\otimes n})$ then the $F$-signature is given by
\[ s(R) = \begin{cases}\frac{6d^2-c^2}{6d^2(ad+c)} & \text{ if $c\leq d$} \\
\frac{6d^3-(2d-c)^3}{6cd^2(ad+c)}&  \text{ if $ d \leq c \leq 2d$}\\
\frac{d}{c(ad+c)} & \text{ if $  2d\leq c$}.
\end{cases}\]
\end{Theoremx}

\section{Preliminaries}\label{sec:preliminaries}
\subsection{Properties of the \texorpdfstring{$F$}{F}-signature function} In this subsection, we recall how the $F$-signature function is defined on the ample cone of a projective variety and summarize some of its properties. Throughout this subsection we let $X$ be a normal projective variety over an algebraically closed field $k$ of positive characteristic $p>0$.

\begin{definition}
A normal projective variety $X$ is said to be \emph{globally $F$-regular} if for any effective Weil divisor $D$ on $X$, there exists an integer $e>0$ such that map $\cO_X \to F^e_*\cO_X(D)$ splits as a map of $\cO_X$-modules.

\end{definition}

Let $\cL$ be an ample divisor on $X$. Then, we have a section ring $S_{\cL}:= \bigoplus_{n\in \mathbb{Z}} H^0(X,\cL^{\otimes n})$ such that $X\cong \mathrm{Proj}(S_{\cL})$. By \cite{SS10}, $X$ is globally $F$-regular if and only if $S_{\cL}$ is strongly $F$-regular. Combined with the equivalence between strong $F$-regularity and positivity of the $F$-signature \cite[Theorem 0.2]{AL03}, we see that $X$ is globally $F$-regular if and only if there is an ample divisor $\cL$ whose section ring has positive $F$-signature. Furthermore, it is proven in \cite{SS10} that for a globally $F$-regular variety $X$, there is an effective divisor $\Delta$ such that $(X,\Delta)$ is log-Fano.

\begin{definition}\label{def:f-sig-1}
Let $X$ be a globally $F$-regular projective variety over $k$ and let $\cL$ be an ample divisor on $X$. Letting $S_{\cL}$ be the section ring of $\cL$ as above, the $F$-signature of $\cL$ is defined to be $s(\cL):= s(S_\cL)$.
\end{definition}
The $F$-signature function is well-defined on the ample classes and is Lipschitz continuous on the ample cone of $X$. Furthermore, the $F$-signature function extends to $\mathrm{Nef}_{\mathbb{R}}(X)\setminus \{0\}$ and is zero if $\cL\in \mathrm{Nef}_\mathbb{R}(X)$ is not big \cite{LP24}. Additionally, the $F$-signature function enjoys the following transformation rule (see \cite{VK12} for Veronese subrings and \cite{CR22} for more general types of extensions). 

\begin{theorem}\cite[Theorem 2.6.2]{VK12}\label{thm:FsigTransRule} Let $X$ be a projective variety over $k$ and $\cL$ be an ample invertible divisor on $X$. Let $S_\cL$ and $S_{\cL^n}$ denote the section rings with respect to $\cL$ and $\cL^n$ respectively for $n\in \mathbb{N}$. Then, 
\[s(S_{\cL}) = n\cdot s(S_{\cL^n}).\]
\end{theorem}

\begin{remark}\cref{thm:FsigTransRule} suggests that if an ample divisor $\cL$ is a finite sum $a_1 e_{i_1}+\dots +a_te_{i_t}$ of basis elements of the N\'eron--Severi group such that each $e_{i_k}$ is an ample class and the $F$-signature function $s(\cL)$ is a rational function on the coefficients of $e_{i_k}$, then $s(\cL)$ must be a homogeneous function of degree $-1$. This demonstrates an issue with the computation of the $F$-signature function on the Hirzebruch surface in \cite{HS17} (the cases on pages 14 to 16 of \emph{op. cit.}, for instance,  do not pass this criterion).
\end{remark}\label{rem:transformation-rule}

\subsection{The ample cone of a Hirzebruch surface}
Let $X = \mathcal{H}_a$ be the Hirzebruch surface of degree $a$ and let $\cL$ be an ample divisor on $X$. Since $X$ is a projective toric variety, the main tool we employ to compute $s(\cL)$ is Von Korff's theorem concerning the volume of the associated polytope. To state Von Korff's result, we first recall some conventions and establish the notation that will be in use throughout the article, primarily following \cite{CLS11,Ful93}. Consider the following figure:
\begin{figure}[H]
\centering
\begin{tikzpicture}
\draw[<->]   (-3,0) -- (3, 0);
\draw[<->]    (0,-3) -- (0,3);

\draw[->,ultra thick] (0,0) -- (0,1) node[anchor=north west] {$\vec{v}_2$};
\draw[->,ultra thick] (0,0) -- (1,0) node[anchor=north] {$\vec{v}_1$};
\draw[->,ultra thick] (0,0) -- (-1,2) node[anchor=south] {$\vec{v}_3$};
\draw[->,ultra thick] (0,0) -- (0,-1) node[anchor=west] {$\vec{v}_4$};

\draw[dotted] (0,2) -- (-1,2);
\draw[dotted] (-1,0) -- (-1,2);

\end{tikzpicture}
\caption{Primitive generators of the Hirzebruch surface.}
\label{FigPrimGenOfHirzebruchSurf}
\end{figure}
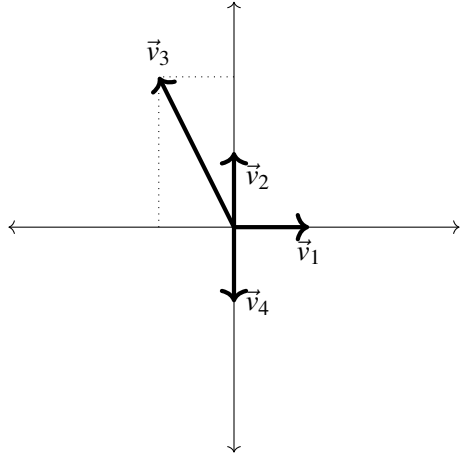
\cref{FigPrimGenOfHirzebruchSurf} describes the primitive generators of the fan of $X$. Let $\vec{v}_1 = \vec{e}_1$, $\vec{v}_2 = \vec{e}_2$, $\vec{v}_3 = -\vec{e}_1+a\vec{e}_2$, and $\vec{v}_4 = -e_2$. Consider the torus-invariant divisors $D_1$, $D_2$, $D_3$, and $D_4$ corresponding to each generator. For a divisor $D = \sum_{i=1}^{4} b_i D_i$, recall that the associated polytope $P_D$ is defined as \[\begin{split} P_D &= \{ \vec{u}\in M_\mathbb{R} \mid \langle \vec{u}, \vec{v}_j\rangle \geq - b_j \text{ for all $j$}\}\\ & =\{ \vec{u}\in \mathbb{R}^2 \mid  \vec{u}\cdot \vec{v}_j \geq - b_j \text{ for all $j$}\}.\end{split}\] Furthermore, consider the polytope $P_{\Sigma, D}$ given by the set of all $(\vec{u}, t)\in \mathbb{R}^2\times \mathbb{R}$ such that $0 \leq (\vec{u},t)\cdot (\vec{v}_j, b_j)\leq 1$ for all $j$ such that $H_j =\{\vec{v}\in\mathbb{R}^2\mid \vec{v}\cdot \vec{v}_j= -b_j\}$ defines a facet of $P_D$. Von Korff's result applied to our setting then says the following:
\begin{theorem}[{\cite[Theorem 4.4.4]{VK12}}]
    With the notation as above, if $D$ is ample then $s_X(D) = \mathrm{Vol}(P_{\Sigma,D})$.
\end{theorem}

Now, since $\div(\chi^{\vec{e}_1})$ and $\div(\chi^{\vec{e}_2})$ are linearly equivalent to $0$, we have \[\begin{split} 0 \sim \div(\chi^{\vec{e}_1}) = \sum_{i=1}^4 \langle e_1, \vec{v}_i \rangle D_i &= D_1 - D_3 \\ 0 \sim \div(\chi^{\vec{e}_2}) = \sum_{i=1}^4 \langle e_2, \vec{v}_i \rangle D_i& = D_2 + aD_3-D_4\end{split}\] Hence, $D_1 \sim D_3$ and $D_2 \sim D_4 - aD_3$. Hence, any divisor $\cL$ on $X$ can be written as \[\begin{split} \cL & = b_1 D_1 + b_2 D_2 + b_3 D_3 + b_4 D_4 \\ & \sim \underbrace{(b_1 -ab_2 + b_3)}_{\colon=c} D_1 + \underbrace{(b_2+b_4)}_{\colon=d}D_4\end{split}\] Here, note that $\cL$ is ample if and only if $c>0$ and $d>0$. Under this notation, the Hilbert--Kunz functions have been computed by Trivedi.
\begin{theorem}\cite[Theorem 3.2]{Tri16}\label{theorem:Trivedi}
     With the notation as above and of \cref{theorem:B}, we have
\begin{align*}
    \ehkx(c,d) = \begin{cases}
        \Omega:=\left(c+\frac{ad}{2}\right)\left(\frac{d}{3}+\frac{(d+1)d}{6c(ad+c)}+\frac{1}{2}+\frac{1}{6d}\right): & c\geq d\\
        \Omega + d\left(c+\frac{(d+1)a}{2}\right)\left(\frac{(a+1)^3}{6a(ad+c)}-\frac{1}{6ac}-\frac{a}{6d}-\frac{1}{2d}+\frac{c}{6d^2}\right): & c< d
    \end{cases}
\end{align*}
\end{theorem}

We turn now and in \cref{sec:calculation} to the calculation of the $F$-signature. Fix $\cL = cD_1 + dD_4$. Specializing \cref{def:f-sig-1} to this setting, we denote $$s_X(c,d):=s(\cL).$$ Here, when $a=1$, $\mathcal{H}_1$ is $\mathbb{P}^2$ blown up at one point. Note that $D_1$ is linearly equivalent to the strict transform of a line in $\mathbb{P}^2$ passing through the blown-up point (see \cite[Theorem 10.4.4]{CLS11}). Further, $D_4$ is linearly equivalent to the exceptional divisor. The $F$-signature function on the ample cone in this setting has been handled already in see \cite{LP24} and \cite[Appendix A]{Lee23}, and is given as follows:

\begin{align*}
    s_X(c,d)=\begin{cases}
        \frac{d}{c(c-d)}: & 3d\leq c\\
        \frac{c-2d}{2cd}+\frac{(2c-3d)(3d-c)}{6d^2(c-d)}+\frac{(3d-c)^2}{2cd^2}: & 2d \leq c \leq 3d\\
        \frac{1}{c}-\frac{(c-d)^3+(2d-c)^3}{6cd^2(c-d)}: & \frac{3}{2}d\leq c\leq 2d\\
        \frac{1}{c}-\frac{(c-d)^2+(2d-c)^2+(3d-2c)(2d-c)+(3d-2c)^2}{6cd^2}: & 0\leq c\leq \frac{3}{2}d
    \end{cases}
\end{align*}

\section{Proof of \texorpdfstring{\cref{theorem:B}}{Theorem B}}\label{sec:calculation}

Continuing with the notation of \cref{sec:preliminaries}, fix for the rest of this section an integer $a\geq 2$. To compute the $F$-signature function of $\mathcal{H}_a$, we compute $P_{\Sigma, L}$ as follows.

\[\begin{split} P_{\Sigma, L} &= \left\{(x,y,z) \middle|\, \begin{matrix}0 \leq (x,y,z)\cdot (\vec{v}_1, c) < 1\\0\leq (x,y,z)\cdot (\vec{v}_2,0) < 1\\0\leq (x,y,z)\cdot (\vec{v}_3,0) < 1\\0\leq (x,y,z)\cdot (\vec{v}_4, d) < 1\end{matrix}\right\}\\ & = \left\{(x,y,z) \middle|\, \begin{split}0\leq x+cz < 1,&\quad  0\leq y < 1,\\ 0\leq -x+ay < 1,& \quad 0\leq -y+zd<1 \end{split}\right\}.\end{split}\]

Now, after the change of coordinates $x \mapsto -x+ay$, $y\mapsto y$, $z\mapsto z$, the volume of $P_{\Sigma, L}$ is the volume of the polytope enclosed in 

\renewcommand{\theequation}{\roman{equation}}
\begin{numcases}{}  0\leq x < 1\label{polytope:i}\\ 0\leq y< 1\label{polytope:ii}\\0\leq -y+zd<1\label{polytope:iii}\\0\leq -x+ay+cz < 1\label{polytope:iv}\end{numcases}

\begin{figure}[H]
\centering
\begin{tikzpicture}[
Mark/.style={ circle, draw, inner sep=0pt, minimum size=4pt, fill=black},
]

\draw[dotted,thick]   (0,0,0) -- (2, 0,0);

\draw[->]   (2,0,0) -- (4, 0,0) node[anchor=north] {$y$};
\draw[->]    (0,1,0) -- (0,3,0) node[anchor=north east] {$z$};

\draw[->](0,0,3) -- (0,0,4.5) node[anchor=north east] {$x$};

\draw[-,ultra thick]   (0,0,3) -- (3, 1,3);

\draw[dotted,ultra thick]   (0,0,0) -- (3, 1,0);

\draw[-,ultra thick]   (3, 1,3) -- (3, 1,0);

\draw[dotted,ultra thick]   (0, 0,0) -- (0, 0,3);

\draw[-,ultra thick]   (0, 0,3) -- (0, 1,3);

\draw[dotted,ultra thick]   (0, 0,0) -- (0, 1,0);

\draw[-,ultra thick]   (0, 1,3) -- (0, 1,0);

\draw[-,ultra thick]   (0,1,3) -- (3, 2,3);

\draw[-,ultra thick]   (0,1,0) -- (3, 2,0);

\draw[-,ultra thick]   (3, 2,3) -- (3, 2,0);

\draw[-,ultra thick]   (0, 1,0) -- (0, 1,3);

\draw[-,ultra thick]   (3, 1,3) -- (3, 2,3);

\draw[-,ultra thick]   (3, 1,0) -- (3, 2,0);

\draw[dotted] (3,3.5,3) -- (3,0,3);

\draw[dotted] (3,0,3) -- (3,0,0);

\draw[dotted] (3,0,3) -- (0,0,3);

\draw[dotted] (0,0,3) -- (0,3.5,3);
\draw[dotted] (3,0,0) -- (3,3.5,0);

\node at (-.2,0,1.5) {\scriptsize\ref{edge:1}};
\node at (1.7,0.2,3) {\scriptsize\ref{edge:2}};
\node at (3,.75,1.5) {\scriptsize\ref{edge:3}};
\node at (1.3,0.2,0) {\scriptsize\ref{edge:4}};

\node at (-0.3,1,1.5) {\scriptsize\ref{edge:5}};
\node at (1.7,1.7,3) {\scriptsize\ref{edge:6}};
\node at (3,1.7,1.5) {\scriptsize\ref{edge:7}};
\node at (1.3,1.7,0) {\scriptsize \ref{edge:8}};

\node at (-0.3,0.5,3) {\scriptsize\ref{edge:9}};
\node at (3,0.7,3) {\scriptsize\ref{edge:10}};
\node at (3.4,1.5,0) {\scriptsize\ref{edge:11}};
\node at (0,1.3,0) {\scriptsize\ref{edge:12}};

\node[Mark] at (3,0,0) {};
\node at (3.3,-0.4,0.4) {$(0,1,0)$};

\node[Mark] at (0,0,3) {};
\node at (0.4,-0.4,3) {$(1,0,0)$};
\end{tikzpicture}
\caption{The graph of the polytope $P$ enclosed in \eqref{polytope:i}, \eqref{polytope:ii}, and \eqref{polytope:iii}, with edge labels.}
\label{GraphofP:labels}
\end{figure}
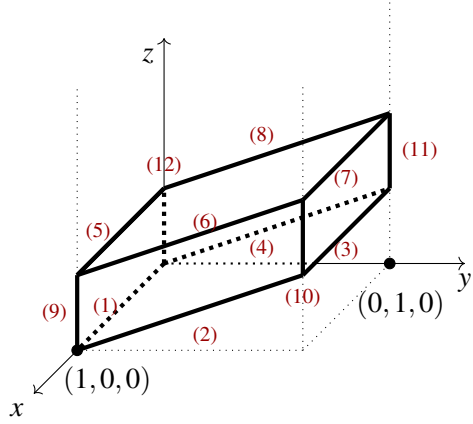

To compute the volume, we fix the range enclosed in \eqref{polytope:i}, \eqref{polytope:ii}, and \eqref{polytope:iii} and find the intersection with \eqref{polytope:iv} in different cases. Let $P$ be the polytope enclosed in \eqref{polytope:i}, \eqref{polytope:ii}, and \eqref{polytope:iii}. The $F$-signature is then the volume of the intersection of $P$ and \eqref{polytope:iv}. To identify the intersection with the range in \eqref{polytope:iv}, we find the intersection of the edges of the polytope $P$ with extremes of \eqref{polytope:iv} (i.e., the intersection of the edges of the polytope $P$ with the two planes $0=-x+ay +cz$ and $1 = -x+ay +cz$). We label \ref{edge:1} - \ref{edge:4} as edges on the bottom of the polytope, \ref{edge:5} - \ref{edge:8} as the ones for the top lid, and \ref{edge:9} - \ref{edge:12} as the sides of the polytope, as depicted in \cref{GraphofP:labels}.

The following table is the intersection of the polytope with the extreme planes of \eqref{polytope:iv}. $\emptyset$ in the table indicates that the intersection does not exist.

\renewcommand{\arraystretch}{1.5}
\begin{table}[H]
\begin{center}
\begin{tabular}{|N|c|c|c| } 
  \hline
  \multicolumn{1}{|c|}{\multirow{2}{*}{Edge ($n$) in \cref{GraphofP:labels}}} & \multirow{1}{*}{Point $v_n$ of intersection between} & \multirow{1}{*}{Point $u_n$ of intersection between}\\
 \multicolumn{1}{|c|}{} & \multirow{1}{*}{edge ($n$) and plane $0 = -x+ay+cz$} & \multirow{1}{*}{edge ($n$) and plane $1 = -x+ay+cz$}\\
  \hline
   \label{edge:1}& $(0,0,0)$ & $\emptyset$\\ 
  \hline
   \label{edge:2}& $\left(1, \frac{d}{ad+c}, \frac{1}{ad+c}\right)$ & $\left(1, \frac{2d}{ad+c}, \frac{2}{ad+c}\right)$ \\ 
  \hline
 \label{edge:3}& $\emptyset$  & $\emptyset$  \\ 
  \hline
    \label{edge:4} & $(0,0,0)$ & $\left(0, \frac{d}{ad+c}, \frac{1}{ad+c}\right)$ \\ 
  \hline
   \label{edge:5} & $\left(\frac{c}{d}, 0, \frac{1}{d}\right)$ & $\left(\frac{c}{d}-1, 0, \frac{1}{d}\right)$ \\ 
  \hline
  \label{edge:6}  & $\left(1, \frac{d-c}{ad+c}, \frac{a+1}{ad+c}\right)$ & $\left(1, \frac{2d-c}{ad+c}, \frac{a+2}{ad+c}\right)$ \\ 
  \hline
   \label{edge:7}& $\emptyset$  & $\emptyset$ \\ 
  \hline
   \label{edge:8}& $\emptyset$ &$\left(0, \frac{d-c}{ad+c}, \frac{a+1}{ad+c}\right)$  \\ 
  \hline
   \label{edge:9} & $\left(1,0, \frac{1}{c}\right)$ & $\left(1,0, \frac{2}{c}\right)$ \\ 
  \hline
    \label{edge:10} & $\emptyset$  &  $\emptyset$\\ 
  \hline
   \label{edge:11} & $\emptyset$ & $\emptyset$ \\ 
  \hline
   \label{edge:12} & $(0,0,0)$ & $\left(0,0,\frac{1}{c}\right)$  \\ 
  \hline
\end{tabular}
\end{center}
\caption{Points of intersection $v_n$ and $u_n$ (if they exist) between edges ($n$) of the polytope and the planes $0 = -x+ay+cz$ and $1 = -x+ay+cz$ respectively.}
\label{IntersectionPointsPolytope}
\end{table}

Now, we analyze which ranges of $c$ and $d$ make the intersection points present on the polytope $P$. Here, note that we cut the square prism with two planes in \cref{GraphofP:labels}. Hence, the $x$ and $y$ coordinates in \cref{IntersectionPointsPolytope} must lie between $0$ and $1$. Furthermore, the $z$-component of the bottom edges ranges from $0$ to $1/d$ and the top edges range from $1/d$ to $2/d$.

For $v_2$ to lie on $P$, we must have $0 \leq \frac{d}{ad+c} \leq 1$. This is always possible since we have $d>0$ and $d \leq ad+c$. Similarly for $u_2$, we must have $0\leq \frac{2d}{ad+c} \leq 1$ which is guaranteed since $2d\leq ad+c$.

For $v_5$, we must have $0\leq \frac{c}{d}\leq 1$. This point is present only when $c\leq d$. Similarly for $u_5$, we need $0\leq \frac{c}{d}-1\leq 1$. Hence, this intersection point exists when $d\leq c \leq 2d$.

For $v_6$, we must have $0\leq \frac{d-c}{ad+c}\leq 1$ and $\frac{1}{d}\leq \frac{a+1}{ad+c}\leq \frac{2}{d}$. The two ranges are valid only when $c\leq d$. On the other hand, $u_6$ requires $0\leq \frac{2d-c}{ad+c}\leq 1$ and $\frac{1}{d}\leq \frac{a+2}{ad+c}\leq \frac{2}{d}$. This point is present when $c\leq 2d$. Furthermore, the point $u_8$ has the same $y$ and $z$ coordinates as $v_6$, hence it is present when $c\leq d$. 

The point $v_9$ exists when $1/c \leq 1/d$, or $d\leq c$ (and so does $u_{12}$). The point $u_9$ exists when $\frac{2}{c} \leq \frac{1}{d}$, or equivalently, $2d\leq c$. Hence, there are three different cases to handle. The table below summarizes the intersections we have to consider for the vertices of $P_{\Sigma,L}$).

\begin{table}[H]
\begin{center}
\begin{tabular}{ |c|c|c|c| } 
 \hline
  & \cref{case:1}: $c\leq d$ & \cref{case:2}: $d\leq c\leq 2d$ & \cref{case:3}: $2d\leq c$ \\
 \hline
 $0=-x+ay+cz$ & $v_1$, $v_2$, $v_5$, $v_6$ & $v_1$, $v_2$, $v_9$ & $v_1$, $v_2$, $v_9$\\ 
 \hline
 $1=-x+ay+cz$ & $u_2$, $u_4$, $u_6$, $u_8$ & $u_2$, $u_4$, $u_5$, $u_6$, $u_{12}$ & $u_2$, $u_4$, $u_9$, $u_{12}$\\
 \hline
\end{tabular}
\end{center}
\caption{Vertices of $P_{\Sigma, L}$}
\label{Vertices}
\end{table}

\begin{case}[$c\leq d$]\label{case:1}
The strategy here is that we will subtract the parts belonging to $ -x+ay+cz\leq 0$ or $1\leq -x+ay+cz $ from the polytope $P$. Consider first the intersection of $P$ with the range $-x+ay+cz\leq 0$. As described in the picture below, it is a triangular prism with six vertices including the intersection as in \cref{Vertices}.

\begin{figure}[H]
\centering
\begin{tikzpicture}[
Mark/.style={ circle, draw=cerulean, inner sep=0pt, minimum size=4pt, fill=cerulean},
]
\draw[dotted,thick]   (0,0,0) -- (2, 0,0);

\draw[->]   (2,0,0) -- (4, 0,0) node[anchor=north] {$y$};
\draw[->]    (0,1,0) -- (0,3,0) node[anchor=north east] {$z$};

\draw[->](0,0,3) -- (0,0,4.5) node[anchor=north east] {$x$};

\draw[-,ultra thick,cerulean]   (0,0,3) -- (3*3/8, 3/8 ,3);
\draw[-,ultra thick]   (3*3/8, 3/8 ,3) -- (3, 1,3);

\draw[dotted,ultra thick]   (0,0,0) -- (3, 1,0);

\draw[-,ultra thick]   (3, 1,3) -- (3, 1,0);

\draw[dotted,ultra thick,cerulean]   (0, 0,0) -- (0, 0,3);

\draw[-,ultra thick,cerulean]   (0, 0,3) -- (0, 1,3);

\draw[dotted,ultra thick]   (0, 0,0) -- (0, 1,0);

\draw[-,ultra thick]   (0, 1,3) -- (0, 1,0);

\draw[-,ultra thick,cerulean]   (0,1,3) -- (3*1/8,3*3/8,3);
\draw[-,ultra thick]   (3*1/8,3*3/8,3) -- (3, 2,3);

\draw[-,ultra thick]   (0,1,0) -- (3, 2,0);

\draw[-,ultra thick]   (3, 2,3) -- (3, 2,0);

\draw[-,ultra thick]   (0, 1,0) -- (0,1,3*2/3);
\draw[-,ultra thick,cerulean]   (0,1,3*2/3) -- (0, 1,3);

\draw[-,ultra thick]   (3, 1,3) -- (3, 2,3);

\draw[-,ultra thick]   (3, 1,0) -- (3, 2,0);

\node[Mark] at (0,0,0) {}; 

\node[Mark] at (3*3/8, 3/8 ,3) {};

\node[Mark] at (0,1,3*2/3) {};

\node[Mark] at (3*1/8,3*3/8,3) {}; 

\node[Mark] at (0,1,3) {}; 

\node[Mark] at (0,0,3) {}; 

\draw[-,cerulean,thick]  (0,0,0) -- (3*3/8, 3/8 ,3); 
\draw[-,cerulean,thick]  (3*3/8, 3/8 ,3) -- (3*1/8,3*3/8,3); 

\draw[-,cerulean,thick]  (3*1/8,3*3/8,3) -- (0,1,3*2/3); 

\draw[-,cerulean,thick]  (0,1,3*2/3)-- (0,0,0);
\end{tikzpicture}\hspace{2cm}
\begin{tikzpicture}[
Mark/.style={ circle, draw=cerulean, inner sep=0pt, minimum size=4pt, fill=cerulean},
]
 
\node[Mark] at (0,0,0) {};
\node at (0+0.4,0,0) {$v_1$};

\node[Mark] at (2*3*3/8, 2*3/8 ,2*3) {}; 
\node at (2*3*3/8+0.4, 2*3/8 ,2*3) {$v_2$};

\node[Mark] at (0, 2*1, 2*3*2/3) {};
\node at (0, 2*1+0.4, 2*3*2/3) {$v_5$};

\node[Mark] at (2*3*1/8, 2*3*3/8, 2*3) {};
\node at (2*3*1/8, 2*3*3/8-0.4, 2*3) {$v_6$};

\node[Mark] at (0,2*1,2*3) {}; 
\node at (0,2*1 ,2*3+0.7) {$B$};

\node[Mark] at (0,0,2*3) {}; 
\node at (0,0,2*3+0.7) {$A$};

\draw[-,ultra thick,cerulean]  (0,0,0) -- (2*3*3/8, 2*3/8 , 2*3); 
\draw[-,ultra thick,cerulean]  (2*3*3/8, 2*3/8 , 2*3) -- (2*3*1/8, 2*3*3/8, 2*3); 

\draw[-,ultra thick,cerulean]  (2*3*1/8, 2*3*3/8, 2*3) -- (0, 2*1, 2*3*2/3); 

\draw[dotted,ultra thick,cerulean]  (0,0,2*3)-- (0,0,0); 

\draw[-,ultra thick,cerulean]  (0,0,2*3)-- (2*3*3/8, 2*3/8 ,2*3); 

\draw[-,ultra thick,cerulean]  (0,2*1,2*3*2/3)-- (0,0,0); 

\draw[-,ultra thick,cerulean] (0,0,2*3) -- (0,2*1,2*3);

\draw[-,ultra thick,cerulean] (0,2*1,2*3) -- (2*3*1/8, 2*3*3/8, 2*3);

\draw[-,ultra thick,cerulean] (0,2*1,2*3) -- (0,2*1, 2*3*2/3);

\end{tikzpicture}
\caption{The intersection of the polytope $P$ with the region $-x+ay+cz\leq 0$ when $c\leq d$.}\label{fig:case-1-prism}
\end{figure}
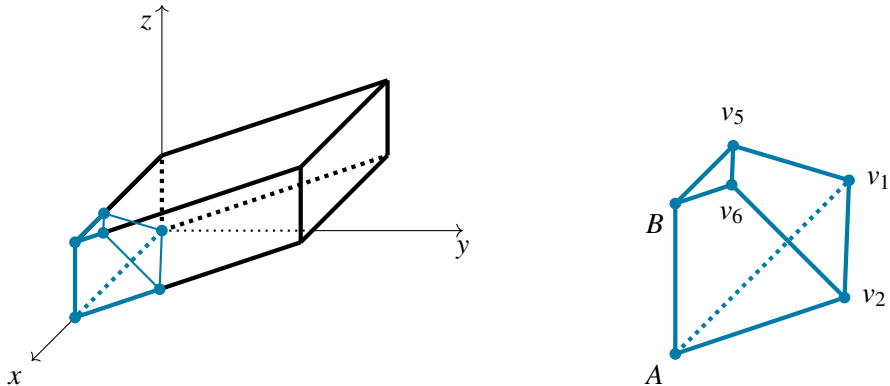

We can triangulate the prism with three triangular pyramids. Let $P_1$ be the pyramid whose vertices are $v_1$, $v_2$, $A$, $B$; $P_2$ be the one with $v_1$, $v_5$, $v_6$, $B$; and $P_3$ be the one with $v_1$, $v_2$, $v_6$, $B$ as in \cref{fig:case-1-prism}. Then, the volume of the prism in \cref{fig:case-1-prism} is $\vol(P_1)+\vol(P_2)+\vol(P_3)$ where \[\begin{split}\vol(P_1)& = \frac{1}{6}\begin{vmatrix} 1&0&0&0 \\ 1&1&0&0 \\ 1&1 & \frac{d}{ad+c} & \frac{1}{ad+c} \\ 1 & 1& 0 & \frac{1}{d}\end{vmatrix} = \frac{1}{6(ad+c)} \\ \vol(P_2) &= \frac{1}{6}\begin{vmatrix} 1&0&0&0 \\ 1 & 1 & \frac{d-c}{ad+c} & \frac{a+1}{ad+c} \\ 1 & \frac{c}{d}& 0 & \frac{1}{d} \\ 1&1&0 & \frac{1}{d} \end{vmatrix}=\frac{(d-c)^2}{6d^2(ad+c)}\\
\vol(P_3) &= \frac{1}{6}\begin{vmatrix}1&0&0&0 \\ 1&1 & \frac{d}{ad+c} & \frac{1}{ad+c} \\ 1&1&\frac{d-c}{ad+c} & \frac{a+1}{ad+c} \\ 1&1&0 & \frac{1}{d} \end{vmatrix} = \frac{d-c}{6d(ad+c)}
\end{split}\]
Hence, 
\begin{align*}\vol(P_1)+\vol(P_2)+\vol(P_3) = \frac{(d-c)^2 + d(2d-c)}{6d^2(ad+c)}.
\end{align*}

Now, consider the intersection of the range $1\leq -x+ay+cz $ with $P$. The polytope of the intersection is as follows:

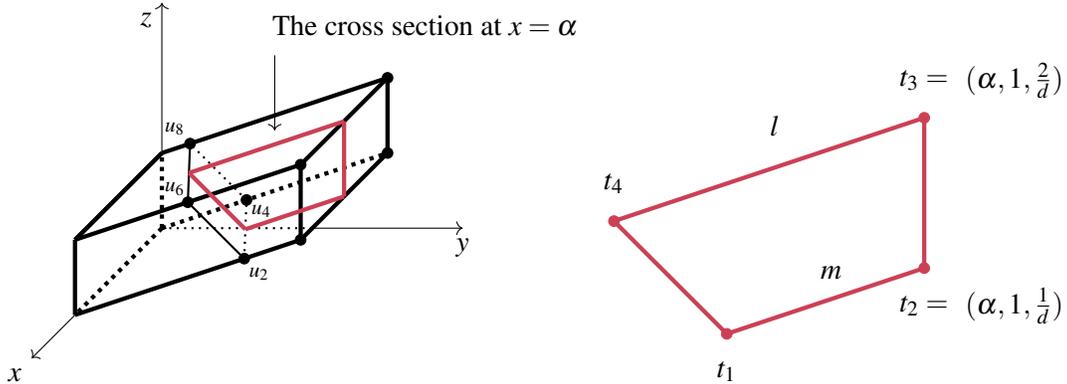
\begin{figure}[H]
\begin{tikzpicture}[
Mark/.style={ circle, draw, inner sep=0pt, minimum size=4pt, fill=black}]

\draw[dotted,thick]   (0,0,0) -- (2, 0,0);

\draw[->]   (2,0,0) -- (4, 0,0) node[anchor=north] {$y$};
\draw[->]    (0,1,0) -- (0,3,0) node[anchor=north east] {$z$};

\draw[->](0,0,3) -- (0,0,4.5) node[anchor=north east] {$x$};

\draw[-,ultra thick]   (0,0,3) -- (3, 1,3);

\draw[dotted,ultra thick]   (0,0,0) -- (3, 1,0);

\draw[-,ultra thick]   (3, 1,3) -- (3, 1,0);

\draw[dotted,ultra thick]   (0, 0,0) -- (0, 0,3);

\draw[-,ultra thick]   (0, 0,3) -- (0, 1,3);

\draw[dotted,ultra thick]   (0, 0,0) -- (0, 1,0);

\draw[-,ultra thick]   (0, 1,3) -- (0, 1,0);

\draw[-,ultra thick]   (0,1,3) -- (3, 2,3);

\draw[-,ultra thick]   (0,1,0) -- (3, 2,0);

\draw[-,ultra thick]   (3, 2,3) -- (3, 2,0);

\draw[-,ultra thick]   (0, 1,0) -- (0, 1,3);

\draw[-,ultra thick]   (3, 1,3) -- (3, 2,3);

\draw[-,ultra thick]   (3, 1,0) -- (3, 2,0);

\node[Mark] at (3*6/8, 3*2/8 ,3) {};
\node at (3*6/8+.2, 3*2/8-.2,3) {\scriptsize$u_2$};

\node[Mark] at (3*3/8, 3/8 ,0) {};
\node at (3*3/8+.2, 3/8-.15 ,0) {\scriptsize$u_4$};

\node[Mark] at (3*4/8 , 3*4/8,3) {}; 
\node at (3*4/8-.25,3*4/8+.1,2.8) {\scriptsize$u_6$};

\node[Mark] at (3*1/8, 3*3/8, 0) {}; 
\node at (3*1/8-.2, 3*3/8+.2, 0) {\scriptsize$u_8$};

\node[Mark] at (3,1,3) {};

\node[Mark] at (3,2,3) {}; 

\node[Mark] at (3,2,0) {}; 

\node[Mark] at (3,1,0) {}; 

\draw[dotted,thick]  (3*6/8, 3*2/8 ,3) -- (3*3/8, 3/8 ,0); 
\draw[dotted,thick]  (3*3/8, 3/8 ,0) -- (3*1/8, 3*3/8, 0);

\draw[-,thick]  (3*1/8, 3*3/8, 0) -- (3*4/8 , 3*4/8,3); 

\draw[-,thick] (3*4/8 , 3*4/8,3)-- (3*6/8, 3*2/8 ,3); 

\draw[-,ultra thick,brickred] (3*4.5/8 , 3*1.5/8,1.5)-- (3,1,1.5); 
\draw[-,ultra thick,brickred] (3,2,1.5)-- (3,1,1.5); 
\draw[-,ultra thick,brickred] (3*2.5/8 , 3*3.5/8,1.5)-- (3,2,1.5); 
\draw[-,ultra thick,brickred] (3*4.5/8 , 3*1.5/8,1.5)-- (3*2.5/8 , 3*3.5/8,1.5); 

\node at (3.5,2.7,0) {The cross section at $x=\alpha$};
\draw[->] (1.5,2.3,0) -- (1.5,1.3,0);

\end{tikzpicture}
\begin{tikzpicture}[
Mark/.style={ circle, draw=brickred, inner sep=0pt, minimum size=4pt, fill=brickred},
]

\node at (3*2.5/8 + 2.5, 3*3.5/8 + 2, 1.5) {$l$};
\node at (3*4.5/8 + 2.5, 3*1.5/8 + 0.8, 1.5) {$m$};

\draw[-,ultra thick,brickred] (2*3*4.5/8 ,2* 3*1.5/8 , 2*1.5)-- (2*3, 2*1 , 2*1.5); 
\draw[-,ultra thick,brickred] (2*3, 2*2 , 2*1.5)-- (2*3, 2*1 , 2*1.5);
\draw[-,ultra thick,brickred] (2*3*2.5/8 ,2* 3*3.5/8 , 2*1.5)-- (2*3, 2*2 , 2*1.5);
\draw[-,ultra thick,brickred] (2*3*4.5/8 ,2* 3*1.5/8 , 2*1.5)-- (2*3*2.5/8 , 2*3*3.5/8 , 2*1.5); 

\node[Mark] at (2*3*4.5/8 ,2* 3*1.5/8 , 2*1.5) {}; 
\node at (2*3*4.5/8 ,2* 3*1.5/8 -0.5 , 2*1.5) {$t_1$};

\node[Mark] at (2*3, 2*1 , 2*1.5) {}; 
\node at (2*3, 2*1 -0.5, 2*1.5) {$t_2=$}; 
\node at (2*3+1.2, 2*1 -0.5, 2*1.5) {$(\alpha, 1, \frac{1}{d})$}; 

\node[Mark] at (2*3, 2*2 , 2*1.5) {}; 
\node at (2*3, 2*2 +0.5, 2*1.5) {$t_3=$}; 
\node at (2*3 +1.2 , 2*2 +0.5, 2*1.5) {$(\alpha, 1, \frac{2}{d})$};

\node[Mark] at (2*3*2.5/8 ,2* 3*3.5/8 , 2*1.5) {}; 
\node at (2*3*2.5/8 ,2* 3*3.5/8 +0.5, 2*1.5) {$t_4$}; 

\end{tikzpicture}

\caption{The intersection of the polytope $P$ with the region $1\leq-x+ay+cz$ when $c\leq d$.}
\end{figure}
Now, we observe that the cross section at $x=\alpha$ is a trapezoid and that $x$ can vary on the interval $[0,1]$. Furthermore, since we chopped the parallelopiped, two sides of the trapezoid are parallel. Let $l$ and $m$ be the lines that form parallel sides as above. Let $h$ be the height of the trapezoid. Note that $(\alpha, 0,0)$ is on the line $m$ since $m$ is the intersection of $x=\alpha$ and $-y+zd=0$. We will find the projection of $(\alpha,0,0)$ onto $l$ and find the distance between the two points. Now, $l$ is the intersection of $x=\alpha$ and $-y+zd =1 $ and the direction vector of the line is $\left\langle 0,1, \frac{1}{d}\right\rangle$. Next note that since we have $z=(1+y)/d$ and $x=\alpha$,  $\left(\alpha, y, \frac{y+1}{d}\right)$ parametrizes points on $l$. If $\left(\alpha, y, \frac{y+1}{d}\right)$ is the projection of $(\alpha,0,0)$ onto $l$, it satisfies \[ \left(\left\langle\alpha, y, \frac{y+1}{d}\right\rangle-\left\langle\alpha,0,0\right\rangle\right) \cdot \left\langle0,1,\frac{1}{d}\right\rangle = 0.\]
This is satisfied when $y=-\frac{1}{d^2+1}$ so the height is $h=\frac{1}{\sqrt{d^2+1}}$. Now we find the lengths of the sides of the top and of the base. Note that $t_1$ is in the intersection of the $-x+ay+cz = 1$  and the plane $-y+zd=0$ at $x= \alpha$. Hence, we find the coordinate of $t_1$ by solving the following equations. \[\begin{cases} x=\alpha \\  0 = -y+zd \\ 1 = -x+ay+cz \end{cases} \Rightarrow \begin{cases} x = \alpha \\ y = \frac{d(1+\alpha)}{ad+c} \\ z= \frac{1+\alpha}{ad+c}\end{cases}\] Hence, the coordinate of $t_1$ is $\left(\alpha, \frac{(1+\alpha)d}{ad+c}, \frac{1+\alpha}{ad+c} \right)$.
Similarly, $t_4$ is on the intersection of the $ -x+ay+cz = 1 $ and the plane $-y+zd=1$ at $x=\alpha$, so we find its coordinates by solving the following system.  \[\begin{cases} x = \alpha \\ 1 = -y+zd \\ 1 = -x+ay+cz  \end{cases} \Rightarrow  \begin{cases} x = \alpha \\ y = \frac{(1+\alpha)d -c}{ad+c}\\ z =\frac{1+a+\alpha}{ad+c} \\ \end{cases}\] We then obtain values for $m$ and $l$ as follows. \[m=\sqrt{\left(1-\frac{(1+\alpha)d}{ad+c}\right)^2+\left(\frac{1}{d}-\frac{1+\alpha}{ad+c}\right)^2}=\frac{\sqrt{d^2+1}((a-\alpha -1)d+c)}{d(ad+c)}\]
\[l = \sqrt{\left(1-\frac{(1+\alpha)d-c}{ad+c}\right)^2+\left(\frac{2}{d}-\frac{1+a+\alpha}{ad+c}\right)^2}=\frac{\sqrt{d^2+1}((a-\alpha-1)d + 2c)}{d(ad+c)}\]
Therefore the area of the trapezoid is \[\begin{split}A = \frac{l+m}{2}h &= \frac{2(a-\alpha-1)d+3c}{2d(ad+c)}.
\end{split}\]
Therefore, the volume of the intersection polytope is \[V=\int_0^1 \frac{2(a-\alpha-1)d+3c}{2d(ad+c)} \mathrm{d}\alpha =\frac{(2a-3)d+3c}{2d(ad+c)}. \]
Lastly, the volume of the base prism $P$ is $1/d$. Therefore, the $F$-signature is \[\begin{split}s_X(c,d) &= \frac{1}{d}-\frac{(d-c)^2 + d(2d-c)}{6d^2(ad+c)}-\frac{(2a-3)d+3c}{2d(ad+c)} \\&= \frac{6d^2-c^2}{6d^2(ad+c)}.\end{split}\]
\end{case}

\begin{case}[$d\leq c\leq 2d$]\label{case:2}
Based on the intersection points in Table \ref{IntersectionPointsPolytope}, the polytope that we need to find the volume of is displayed below in \cref{case-2-figure}.

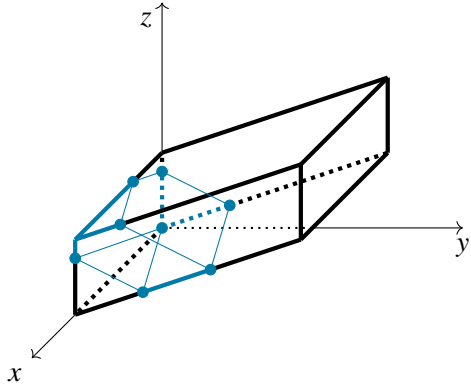
\begin{figure}[H]
\centering
\begin{tikzpicture}[
Mark/.style={ circle, draw=cerulean, inner sep=0pt, minimum size=4pt, fill=cerulean},
]

\draw[dotted,thick]   (0,0,0) -- (2, 0,0);

\draw[->]   (2,0,0) -- (4, 0,0) node[anchor=north] {$y$};
\draw[->]    (0,1,0) -- (0,3,0) node[anchor=north east] {$z$};

\draw[->](0,0,3) -- (0,0,4.5) node[anchor=north east] {$x$};

\draw[-,ultra thick]   (0,0,3) -- (3*3/10, 3*1/10 , 3*1);
\draw[-,ultra thick, cerulean] (3*3/10, 3*1/10 , 3*1) -- (3*6/10, 3*2/10, 3*1);
\draw[-,ultra thick] (3*6/10, 3*2/10, 3*1) -- (3, 1,3);

\draw[dotted,ultra thick,cerulean]   (0,0,0) -- (3*3/10, 3*1/10 ,0);
\draw[dotted,ultra thick] (3*3/10, 3*1/10 ,0) -- (3, 1,0);

\draw[-,ultra thick]   (3, 1,3) -- (3, 1,0);

\draw[dotted,ultra thick]   (0, 0,0) -- (0, 0,3);

\draw[-,ultra thick]   (0, 0,3) -- (0, 3*1/4, 3*1);
\draw[-,ultra thick, cerulean] (0, 3*1/4, 3*1) -- (0, 1,3);

\draw[dotted,ultra thick,cerulean]   (0, 0,0) -- (3*0,3*1/4,3*0);
\draw[dotted,ultra thick] (3*0,3*1/4,3*0) -- (0, 1,0);

\draw[-,ultra thick, cerulean]   (0,1,3) -- (3*2/10,3*4/10,3*1);
\draw[-,ultra thick]  (3*2/10,3*4/10,3*1) -- (3, 2,3);

\draw[-,ultra thick]   (0,1,0) -- (3, 2,0);

\draw[-,ultra thick]   (3, 2,3) -- (3, 2,0);

\draw[-,ultra thick]   (0, 1,0) -- (0,3*1/3,3*1/3);
\draw[-,ultra thick,cerulean]  (0,3*1/3,3*1/3) -- (0, 1,3);

\draw[-,ultra thick]   (3, 1,3) -- (3, 2,3);

\draw[-,ultra thick]   (3, 1,0) -- (3, 2,0);

\coordinate (1) at (0,0,0);
\node[Mark] at (1) {};

\coordinate (2d) at (3*3/10, 3*1/10 , 3*1);
\node[Mark] at (2d) {};

\coordinate (9d) at (0, 3*1/4, 3*1);
\node[Mark] at (9d) {};
\draw[-,cerulean] (1) -- (2d);
\draw[-,cerulean] (2d) -- (9d);
\draw[-,cerulean] (9d) -- (1);

\coordinate (2) at (3*6/10, 3*2/10, 3*1);
\node[Mark] at (2) {};
\coordinate (4) at (3*3/10, 3*1/10 ,0);
\node[Mark] at (4) {};

\coordinate (5) at (0,3*1/3,3*1/3);
\node[Mark] at (5) {};

\coordinate (6) at (3*2/10,3*4/10,3*1);
\node[Mark] at (6) {};

\coordinate (8) at (3*0,3*1/4,3*0);
\node[Mark] at (8) {};

\draw[-,cerulean]  (2) -- (4); 
\draw[-,cerulean]  (4) -- (8); 

\draw[-,cerulean]  (8) -- (5); 

\draw[-,cerulean]  (5)-- (6); 

\draw[-,cerulean] (6) -- (2); 
\end{tikzpicture}
\caption{The polytope for $d\leq c\leq 2d$}\label{case-2-figure}
\end{figure}
For this case, we will refer to the ``top face'' as the face of the polytope that is part of the plane $1=-x+ay+cz$ and the ``bottom face'' as that which is part of $0=-x+ay+cz$. We may easily compute the volume below the bottom face since it is a tetrahedron. Hence, to find the volume of the polytope, we will first find the volume of the polytope enclosed in \eqref{polytope:i}, \eqref{polytope:ii}, \eqref{polytope:iii}, and $-x+ay+cz\leq 1$. Then, we subtract the volume of the polytope enclosed in \eqref{polytope:i}, \eqref{polytope:ii}, \eqref{polytope:iii}, and $-x+ay+ca \leq 0$.
\cref{case2:bottom-face} depicts the polytope without the bottom triangular face:

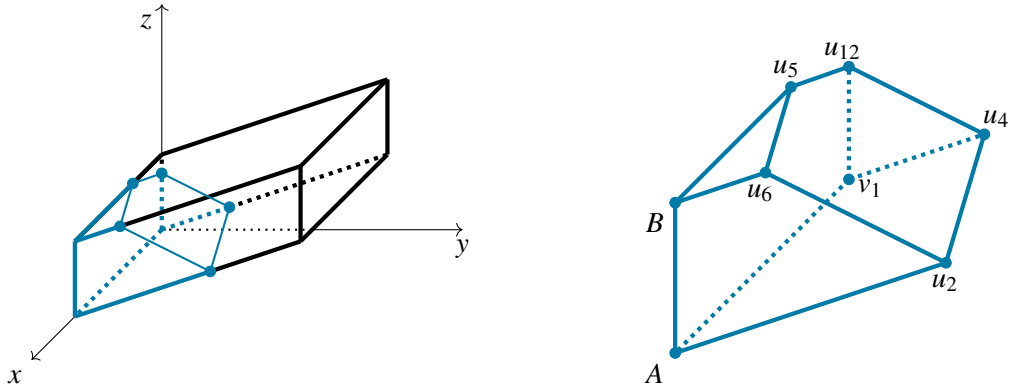
\begin{figure}[H]
\centering
\begin{tikzpicture}[
Mark/.style={ circle, draw=cerulean, inner sep=0pt, minimum size=4pt, fill=cerulean},
]

\draw[dotted,thick]   (0,0,0) -- (2, 0,0);

\draw[->]   (2,0,0) -- (4, 0,0) node[anchor=north] {$y$};
\draw[->]    (0,1,0) -- (0,3,0) node[anchor=north east] {$z$};

\draw[->](0,0,3) -- (0,0,4.5) node[anchor=north east] {$x$};

\draw[-,ultra thick,cerulean]   (0,0,3) -- (3*6/10, 3*2/10, 3*1);
\draw[-,ultra thick]   (3*6/10, 3*2/10, 3*1) -- (3, 1,3);

\draw[dotted,ultra thick,cerulean]   (0,0,0) -- (3*3/10, 3*1/10 ,0);
\draw[dotted,ultra thick] (3*3/10, 3*1/10 ,0) -- (3, 1,0);

\draw[-,ultra thick]   (3, 1,3) -- (3, 1,0);

\draw[dotted,ultra thick,cerulean]   (0, 0,0) -- (0, 0,3);

\draw[-,ultra thick,cerulean]   (0, 0,3) -- (0, 1,3);

\draw[dotted,ultra thick,cerulean]   (0, 0,0) -- (3*0,3*1/4,3*0);
\draw[dotted,ultra thick] (3*0,3*1/4,3*0) -- (0, 1,0);

\draw[-,ultra thick]   (0, 1,3) -- (0, 1,0);

\draw[-,ultra thick,cerulean]   (0,1,3) -- (3*2/10,3*4/10,3*1);
\draw[-,ultra thick] (3*2/10,3*4/10,3*1) -- (3, 2,3);

\draw[-,ultra thick]   (0,1,0) -- (3, 2,0);

\draw[-,ultra thick]   (3, 2,3) -- (3, 2,0);

\draw[-,ultra thick]   (0, 1,0) -- (0,3*1/3,3*1/3);
\draw[-,ultra thick,cerulean]   (0,3*1/3,3*1/3) -- (0, 1,3);

\draw[-,ultra thick]   (3, 1,3) -- (3, 2,3);

\draw[-,ultra thick]   (3, 1,0) -- (3, 2,0);

\coordinate (1) at (0,0,0);

\coordinate (2d) at (3*3/10, 3*1/10 , 3*1);

\coordinate (9d) at (0, 3*1/4, 3*1);

\coordinate (2) at (3*6/10, 3*2/10, 3*1);
\node[Mark] at (2) {}; 
\coordinate (4) at (3*3/10, 3*1/10 ,0);
\node[Mark] at (4) {};

\coordinate (5) at (0,3*1/3,3*1/3);
\node[Mark] at (5) {};

\coordinate (6) at (3*2/10,3*4/10,3*1);
\node[Mark] at (6) {}; 

\coordinate (8) at (3*0,3*1/4,3*0);
\node[Mark] at (8) {};

\draw[-,cerulean,thick]  (2) -- (4); 
\draw[-,cerulean,thick]  (4) -- (8); 

\draw[-,cerulean,thick]  (8) -- (5); 

\draw[-,cerulean,thick]  (5)-- (6);

\draw[-,cerulean,thick] (6) -- (2);
\end{tikzpicture}\hspace{2cm}
\begin{tikzpicture}[
Mark/.style={ circle, draw=cerulean, inner sep=0pt, minimum size=4pt, fill=cerulean},
]

\coordinate (1) at (0,0,0);

\coordinate (2d) at (2*3*3/10, 2*3*1/10 , 2*3*1);

\coordinate (9d) at (0, 2*3*1/4, 2*3*1);

\coordinate (2) at (2*3*6/10, 2*3*2/10, 2*3*1);
\node[Mark] at (2) {}; 
\node at (2*3*6/10-.05, 2*3*2/10-.3, 2*3*1-.1) {$u_2$};

\coordinate (4) at (2*3*3/10,2* 3*1/10 ,0);
\node[Mark] at (4) {};
\node at (2*3*3/10+.2,2* 3*1/10+.25,0+.1) {$u_4$};

\coordinate (5) at (0,2*3*1/3,2*3*1/3);
\node[Mark] at (5) {};
\node at (-.1,2*3*1/3+.2,2*3*1/3-.1) {$u_5$};

\coordinate (6) at (2*3*2/10,2*3*4/10,2*3*1);
\node[Mark] at (6) {};
\node at (2*3*2/10+.2,2*3*4/10,2*3*1+.7) {$u_6$};

\coordinate (12) at (2*3*0,2*3*1/4,2*3*0);
\node[Mark] at (12) {};
\node at (2*3*0-.2,2*3*1/4+.1,2*3*0-.3) {$u_{12}$};

\draw[-,ultra thick, cerulean]  (2) -- (4); 
\draw[-,ultra thick, cerulean]  (4) -- (12); 

\draw[-,ultra thick, cerulean]  (12) -- (5);

\draw[-,ultra thick, cerulean]  (5)-- (6);

\draw[-,ultra thick, cerulean] (6) -- (2);

\coordinate (P1) at (0,0,0);
\node[Mark] at (P1) {};
\node at (0+0.4,0,.3) {$v_1$};

\coordinate (P2) at (0,0,2*3);
\node[Mark] at (P2) {};
\node at (0,0,2*3+0.7) {$A$};

\coordinate (P3) at (0,2*1,2*3);
\node[Mark] at (P3) {};
\node at (0,2*1 ,2*3+0.7) {$B$};

\draw[dotted,ultra thick, cerulean] (P1) -- (P2);
\draw[-,ultra thick, cerulean] (P2) -- (P3);
\draw[dotted,ultra thick, cerulean] (P1) -- (4);
\draw[dotted, ultra thick, cerulean] (P1) -- (12);
\draw[-,ultra thick, cerulean] (P2) -- (2);
\draw[-,ultra thick, cerulean] (P3) -- (5);
\draw[-,ultra thick, cerulean] (P3) -- (6);

\end{tikzpicture}
\caption{The intersection of the polytope $P$ with the region $-x+ay+cz\leq 1$ when $d\leq c\leq 2d$.}
\label{case2:bottom-face}
\end{figure}
Note that the $x$-axis is an edge of the polytope. Hence, we can find the volume of the polytope by integrating the area of the cross-sections parametrized by $x=\alpha$ as in \cref{case:2-triangle,case:2-trapezoid}. The cross-sections are triangular for values $0 \leq x \leq \frac{c}{d}-1$, and trapezoidal for $\frac{c}{d}-1 \leq x \leq 1$. Let us first consider the case when $0 \leq x \leq \frac{c}{d}-1$. Here, we compute the area of the triangle treating the side on the $xz$-plane as a base.
\begin{figure}[H]
\centering
\begin{tikzpicture}[
Mark/.style={ circle, draw, inner sep=0pt, minimum size=4pt, fill=black},
Mark2/.style={ circle, draw=brickred, inner sep=0pt, minimum size=3pt, fill=brickred}
]

\coordinate (1) at (0,0,0);

\coordinate (2d) at (2*3*3/10, 2*3*1/10 , 2*3*1);

\coordinate (9d) at (0, 2*3*1/4, 2*3*1);

\coordinate (2) at (2*3*6/10, 2*3*2/10, 2*3*1);
\node[Mark] at (2) {}; 
\coordinate (4) at (2*3*3/10,2* 3*1/10 ,0);
\node[Mark] at (4) {}; 

\coordinate (5) at (0,2*3*1/3,2*3*1/3);
\node[Mark] at (5) {};

\coordinate (6) at (2*3*2/10,2*3*4/10,2*3*1);
\node[Mark] at (6) {};

\coordinate (12) at (2*3*0,2*3*1/4,2*3*0);
\node[Mark] at (12) {};

\draw[-]  (2) -- (4);
\draw[-]  (4) -- (12); 

\draw[-]  (12) -- (5);

\draw[-]  (5)-- (6);
\draw[-] (6) -- (2);

\coordinate (P1) at (0,0,0);
\node[Mark] at (P1) {};
\coordinate (P2) at (0,0,2*3);
\node[Mark] at (P2) {};
\coordinate (P3) at (0,2*1,2*3);
\node[Mark] at (P3) {};

\draw[dotted,thick] (P1) -- (P2);
\draw[-] (P2) -- (P3);
\draw[dotted,thick] (P1) -- (4);
\draw[dotted, thick] (P1) -- (12);
\draw[-] (P2) -- (2);
\draw[-] (P3) -- (5);
\draw[-] (P3) -- (6);

\coordinate (X) at (0,0,2*4);
\coordinate (XX) at (0,0,-1);
\coordinate (Z) at (0, 2.5, 0);

\draw[->] (P2) -- (X) node[anchor=north west] {$x$};
\draw[dotted, thick] (XX) -- (P1);

\draw[->] (12) -- (Z) node[anchor=south] {$z$};
\coordinate (T1) at (0,0,2*3*1/6);
\coordinate (T2) at (0, 2*3 *7/24 ,2* 3 * 1/6);
\coordinate (T3) at (2*3* 7/20 , 2*3*7/60, 2*3*1/6);
\draw[-, ultra thick, brickred] (T1) -- (T2);
\draw[-, ultra thick, brickred] (T2) -- (T3);
\draw[-, ultra thick, brickred] (T3) -- (T1);

\node[Mark2] at (T1) {};
\node[anchor=north] at (T1) {$(\alpha, 0,0)$};
\node[Mark2] at (T2) {};
\node[anchor=south] at (T2) {$t_2$};
\node[Mark2] at (T3) {};
\node[anchor=west] at (T3) {$t_3$};

\end{tikzpicture}
\caption{Cross sections of the volume in \cref{case2:bottom-face} on the $x$-interval $[0,\frac{c}{d}-1]$.}\label{case:2-triangle}
\end{figure}
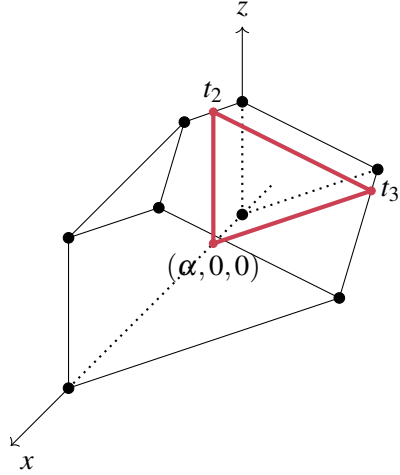
Now, note that $t_2$ is in the intersection of the top face and $xz$-plane. Hence, we find the coordinate of $t_2$ by solving the following equations.
\[\begin{cases} x=\alpha \\ y = 0 \\ 1 = -x+ay+cz \end{cases} \Rightarrow z = \frac{1+\alpha}{c}\]
Similarly, $t_3$ is on the intersection of the top face and the plane $0=-y+zd$ at $x=\alpha$. Hence \[\begin{cases} x = \alpha \\ 1 = -x+ay+cz \\ 0 = -y+zd \end{cases} \Rightarrow \begin{cases} x = \alpha \\ y = \frac{d(1+\alpha)}{ad+c} \\ z= \frac{1+\alpha}{ad+c}\end{cases}\]
Therefore, the area of the cross section is $\frac{1}{2}\frac{1+\alpha}{c}\cdot \frac{d(1+\alpha)}{ad+c}$.

For $\frac{c}{d}-1\leq x\leq 1$, the cross section is a trapezoid. We label $w_1$ and $w_2$ as the vertices not on the $xz$-plane as in \cref{case:2-trapezoid} below.
\begin{figure}[H]
\centering
\begin{tikzpicture}[
Mark/.style={ circle, draw, inner sep=0pt, minimum size=4pt, fill=black},
Mark2/.style={ circle, draw=brickred, inner sep=0pt, minimum size=4pt, fill=brickred}
]

\coordinate (1) at (0,0,0);

\coordinate (2d) at (2*3*3/10, 2*3*1/10 , 2*3*1);

\coordinate (9d) at (0, 2*3*1/4, 2*3*1);

\coordinate (2) at (2*3*6/10, 2*3*2/10, 2*3*1);
\node[Mark] at (2) {}; 
\coordinate (4) at (2*3*3/10,2* 3*1/10 ,0);
\node[Mark] at (4) {}; 

\coordinate (5) at (0,2*3*1/3,2*3*1/3);
\node[Mark] at (5) {}; 

\coordinate (6) at (2*3*2/10,2*3*4/10,2*3*1);
\node[Mark] at (6) {};

\coordinate (12) at (2*3*0,2*3*1/4,2*3*0);
\node[Mark] at (12) {};

\draw[-]  (2) -- (4);
\draw[-]  (4) -- (12);

\draw[-]  (12) -- (5);

\draw[-]  (5)-- (6);

\draw[-] (6) -- (2);

\coordinate (P1) at (0,0,0);
\node[Mark] at (P1) {};
\coordinate (P2) at (0,0,2*3);
\node[Mark] at (P2) {};
\coordinate (P3) at (0,2*1,2*3);
\node[Mark] at (P3) {};

\draw[dotted,thick] (P1) -- (P2);
\draw[-] (P2) -- (P3);
\draw[dotted,thick] (P1) -- (4);
\draw[dotted, thick] (P1) -- (12);
\draw[-] (P2) -- (2);
\draw[-] (P3) -- (5);
\draw[-] (P3) -- (6);

\coordinate (X) at (0,0,2*4);
\coordinate (XX) at (0,0,-1);
\coordinate (Z) at (0, 2.5, 0);

\draw[->] (P2) -- (X) node[anchor=north west] {$x$};
\draw[dotted, thick] (XX) -- (P1);

\draw[->] (12) -- (Z) node[anchor=south] {$z$};

\coordinate (T1) at (2*3* 1/10 , 2*3*11/30, 2*3*2/3);
\coordinate (T2) at (2*3* 5/10 , 2*3*5/30, 2*3*2/3);
\coordinate (T3) at (0,0,2*3*2/3);
\coordinate (T4) at (0, 2*3 *1/3 ,2* 3 * 2/3);

\draw[-, ultra thick, brickred] (T1) -- (T2);
\draw[-, ultra thick, brickred] (T2) -- (T3);
\draw[-, ultra thick, brickred] (T3) -- (T4);
\draw[-, ultra thick, brickred] (T1) -- (T4);

\node at (2*3* 1/10+.4, 2*3*11/30+.05, 2*3*2/3+.05) {$w_1$};
\node[Mark2] at (T1) {};
\node[anchor=west] at (T2) {$w_2$};
\node[Mark2] at (T2) {};
\node[anchor=north] at (T3) {$(\alpha,0,0)$};
\node[Mark2] at (T3) {};
\node[anchor=south east] at (T4) {$(\alpha,0,\frac{1}{d})$};
\node[Mark2] at (T4) {};
\end{tikzpicture}
\caption{Cross sections of the volume in \cref{case2:bottom-face} on the $x$-interval $[\frac{c}{d}-1,1]$.}\label{case:2-trapezoid}
\end{figure}
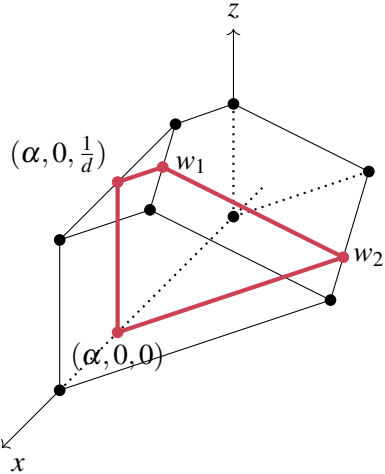

Here, we observe that $l_1$ and $l_2$ are parallel and the distance between the two lines is $\frac{1}{\sqrt{d^2+1}}$ as we have found in \cref{case:1}. Hence, once we find the lengths of $l_1$ and $l_2$, the area of the trapezoid is $\frac{l_1 + l_2}{2\sqrt{d^2+1}}$. In order to find the lengths, we will identify the coordinates of the vertex of the trapezoid.
The only new point of which we need the coordinates is $w_1$. Note that $w_1$ is on the intersection of the top face and $1=-y+zd$. Hence, the coordinate of this point is \[\begin{cases} x=\alpha \\ 1 = -x+ay + cz \\ 1=-y+zd \end{cases} \Rightarrow
 \begin{cases} x = \alpha \\ y = \frac{(1+\alpha)d -c}{ad+c}\\ z =\frac{1+a+\alpha}{ad+c} \\ \end{cases}\] 

Now, $l_1$ is the distance between $(\alpha, 0, \frac{1}{d})$ and $\left(\alpha,  \frac{(1+\alpha)d -c}{ad+c}, \frac{1+a+\alpha}{ad+c} \right)$, which is \[l_1 = \sqrt{\left(\frac{(1+\alpha )d-c}{ad+c}\right)^2+ \left(\frac{1+a+\alpha}{ad+c} - \frac{1}{d}\right)^2}=\sqrt{d^2+1}\cdot\frac{(1+\alpha)d-c}{d(ad+c)}.\]
Similarly, $l_2$ is the distance between $(\alpha, 0,0)$ and $\left(\alpha, \frac{d(1+\alpha)}{ad+c}, \frac{1+\alpha}{ad+c}\right)$: \[l_2 = \sqrt{\frac{(1+\alpha)^2d^2}{(ad+c)^2}+\frac{(1+\alpha)^2}{(ad+c)^2}}=\sqrt{d^2+1}\cdot \frac{1+\alpha}{ad+c}.\]
Therefore, the area of the trapezoid is \[\frac{l_1+l_2}{2\sqrt{d^2+1}} = \frac{2(1+\alpha)d-c}{2d(ad+c)}.\]
Hence, combining the two, we compute the volume of the polytope below $1=-x+ay+cz$ as follows: \[V = \int_0^{c/d-1} \frac{d(1+\alpha)^2}{2c(ad+c)} \mathrm{d}\alpha+ \int_{c/d-1}^1\frac{2(1+\alpha)d-c}{2d(ad+c)} \mathrm{d}\alpha = \frac{c^3 - 6c^2d + 12cd^2 - d^3}{6cd^2(ad+c)}. \]
The polytope below the lower lid is a triangular pyramid made out of the four points \ref{edge:1}, \ref{edge:2}, \ref{edge:9}, and $(1,0,0)$, the volume of which is \[\frac{1}{6}\cdot\begin{vmatrix}1 & 0 & 0 & 0 \\ 1 & 1 & 0 & 0 \\ 1 & 1 & \frac{d}{ad+c} & \frac{1}{ad+c} \\ 1 & 1 & 0 & \frac{1}{c}\end{vmatrix} = \frac{d}{6c(ad+c)}.\]
Hence, the $F$-signature function $s_X(c,d)$ is \[\begin{split}s_X(c,d) &= \frac{c^3 - 6c^2d + 12cd^2 - d^3}{6cd^2(ad+c)} - \frac{d}{6c(ad+c)} \\ &=\frac{6d^3-(2d-c)^3}{6cd^2(ad+c)}. \end{split}\]

\end{case}

\begin{case}[$2d\leq c$]\label{case:3}
For this case, we compute the volume of the polytope displayed in \cref{case-3-fig} below.

\begin{figure}[H]
\centering
\begin{tikzpicture}[
Mark/.style={ circle, draw=cerulean, inner sep=0pt, minimum size=4pt, fill=cerulean},
]

\draw[dotted,thick]   (0,0,0) -- (2, 0,0);

\draw[->]   (2,0,0) -- (4, 0,0) node[anchor=north] {$y$};
\draw[->]    (0,1,0) -- (0,3,0) node[anchor=north east] {$z$};

\draw[->](0,0,3) -- (0,0,4.5) node[anchor=north east] {$x$};

\draw[-,ultra thick]   (0,0,3) -- (3*3/14, 3*1/14 , 3*1);
\draw[-,ultra thick,cerulean]   (3*3/14, 3*1/14 , 3*1) -- (3*6/14, 3*2/14, 3*1);
\draw[-,ultra thick]   (3*6/14, 3*2/14, 3*1) -- (3, 1,3);

\draw[dotted,ultra thick,cerulean]   (0,0,0) -- (3*3/14, 3*1/14 ,0);
\draw[dotted,ultra thick] (3*3/14, 3*1/14 ,0) -- (3, 1,0);

\draw[-,ultra thick]   (3, 1,3) -- (3, 1,0);

\draw[dotted,ultra thick]   (0, 0,0) -- (0, 0,3);

\draw[-,ultra thick]   (0, 0,3) -- (0, 3*1/8, 3*1);
\draw[-,ultra thick,cerulean]   (0, 3*1/8, 3*1) -- (0,3*2/8,3*1);
\draw[-,ultra thick]   (0,3*2/8,3*1) -- (0, 1,3);

\draw[dotted,ultra thick,cerulean]   (0, 0,0) -- (3*0,3*1/8,3*0);
\draw[dotted, ultra thick] (3*0,3*1/8,3*0) -- (0, 1,0);

\draw[-,ultra thick]   (0, 1,3) -- (0, 1,0);

\draw[-,ultra thick]   (0,1,3) -- (3, 2,3);

\draw[-,ultra thick]   (0,1,0) -- (3, 2,0);

\draw[-,ultra thick]   (3, 2,3) -- (3, 2,0);

\draw[-,ultra thick]   (0, 1,0) -- (0, 1,3);

\draw[-,ultra thick]   (3, 1,3) -- (3, 2,3);

\draw[-,ultra thick]   (3, 1,0) -- (3, 2,0);

\coordinate (1) at (0,0,0);
\node[Mark] at (1) {};

\coordinate (2d) at (3*3/14, 3*1/14 , 3*1);
\node[Mark] at (2d) {};

\coordinate (9d) at (0, 3*1/8, 3*1);
\node[Mark] at (9d) {};

\draw[-,cerulean,thick] (1) -- (2d);
\draw[-,cerulean,thick] (2d) -- (9d);
\draw[-,cerulean,thick] (9d) -- (1);

\coordinate (2) at (3*6/14, 3*2/14, 3*1);
\node[Mark] at (2) {};

\coordinate (4) at (3*3/14, 3*1/14 ,0);
\node[Mark] at (4) {};

\coordinate (9) at (0,3*2/8,3*1);
\node[Mark] at (9) {};

\coordinate (12) at (3*0,3*1/8,3*0);
\node[Mark] at (12) {};

\draw[-,cerulean,thick]  (2) -- (4); 
\draw[-,cerulean,thick]  (4) -- (12); 

\draw[-,cerulean,thick]  (9) -- (12);
\draw[-,cerulean,thick]  (9)-- (2); 
\end{tikzpicture}
\caption{The polytope for $2d\leq c$}\label{case-3-fig}
\end{figure}
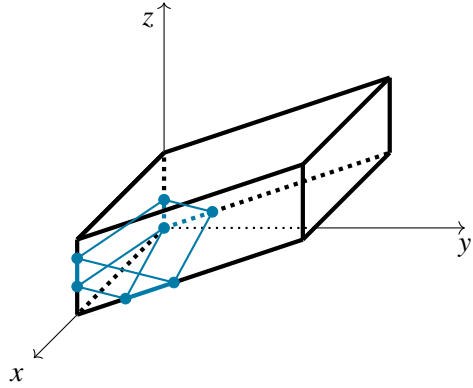
Again, the portion below the bottom face is a tetrahedron with volume given by $\frac{d}{6c(ad+c)}$. We subtract this volume from that of the portion below the top face which, as in \cref{case:2-triangle} of \cref{case:2}, has triangular cross sections.

\begin{figure}[H]
\centering
\begin{tikzpicture}[
Mark/.style={ circle, draw, inner sep=0pt, minimum size=4pt, fill=black},
Mark2/.style={ circle, draw=cerulean, inner sep=0pt, minimum size=4pt, fill=cerulean}
]

\draw[dotted,thick]   (0,0,0) -- (2, 0,0);

\draw[->]   (2,0,0) -- (4, 0,0) node[anchor=north] {$y$};
\draw[->]    (0,1,0) -- (0,3,0) node[anchor=north east] {$z$};

\draw[->](0,0,3) -- (0,0,4.5) node[anchor=north east] {$x$};

\draw[-,ultra thick,cerulean]   (0,0,3) -- (3*6/14, 3*2/14, 3*1);
\draw[-,ultra thick] (3*6/14, 3*2/14, 3*1) -- (3, 1,3);

\draw[dotted,ultra thick,cerulean]   (0,0,0) -- (3*3/14, 3*1/14 ,0);
\draw[dotted,ultra thick]   (3*3/14, 3*1/14 ,0) -- (3, 1,0);
\draw[-,ultra thick]   (3, 1,3) -- (3, 1,0);

\draw[dotted,ultra thick,cerulean]   (0, 0,0) -- (0, 0,3);

\draw[-,ultra thick,cerulean]   (0, 0,3) -- (0,3*2/8,3*1) (0, 1,3);
\draw[-,ultra thick] (0,3*2/8,3*1) -- (0, 1,3);
\draw[dotted,ultra thick,cerulean]   (0, 0,0) -- (3*0,3*1/8,3*0);
\draw[dotted,ultra thick] (3*0,3*1/8,3*0) -- (0, 1,0);

\draw[-,ultra thick]   (0,1,3) -- (3, 2,3);

\draw[-,ultra thick]   (0,1,0) -- (3, 2,0);

\draw[-,ultra thick]   (3, 2,3) -- (3, 2,0);

\draw[-,ultra thick]   (0, 1,0) -- (0, 1,3);

\draw[-,ultra thick]   (3, 1,3) -- (3, 2,3);

\draw[-,ultra thick]   (3, 1,0) -- (3, 2,0);

\coordinate (1) at (0,0,0);

\coordinate (2d) at (3*3/14, 3*1/14 , 3*1);

\coordinate (9d) at (0, 3*1/8, 3*1);

\coordinate (2) at (3*6/14, 3*2/14, 3*1);
\node[Mark2] at (2) {}; 
\coordinate (4) at (3*3/14, 3*1/14 ,0);
\node[Mark2] at (4) {}; 

\node[Mark2] at (0,0,0) {};

\node[Mark2] at (0,0,3) {};

\coordinate (9) at (0,3*2/8,3*1);
\node[Mark2] at (9) {};

\coordinate (12) at (3*0,3*1/8,3*0);
\node[Mark2] at (12) {}; 

\draw[-,cerulean,thick]  (2) -- (4);
\draw[-,cerulean,thick]  (4) -- (12); 

\draw[-,cerulean,thick]  (9) -- (12);
\draw[-,cerulean,thick]  (9)-- (2);
\end{tikzpicture} \hspace{2cm}
\begin{tikzpicture}[
Mark/.style={ circle, draw, inner sep=0pt, minimum size=4pt, fill=black},
Mark2/.style={ circle, draw=cerulean, inner sep=0pt, minimum size=4pt, fill=cerulean},
Mark3/.style={ circle, draw=brickred, inner sep=0pt, minimum size=4pt, fill=brickred}
]

\coordinate (1) at (0,0,0);

\coordinate (2d) at (2*3*3/10, 2*3*1/10 , 2*3*1);

\coordinate (9d) at (0, 2*3*1/4, 2*3*1);

\coordinate (2) at (2*3*6/14, 2*3*2/14, 2*3*1);
\node[Mark2] at (2) {}; 
\coordinate (4) at (2*3*3/14, 2*3*1/14 ,2*0);
\node[Mark2] at (4) {}; 

\coordinate (9) at (2*0,2*3*2/8,2*3*1);
\node[Mark2] at (9) {}; 

\coordinate (12) at (2*3*0,2*3*1/8,2*3*0);
\node[Mark2] at (12) {}; 

\draw[-,cerulean,ultra thick]  (2) -- (4); 
\draw[-,cerulean,ultra thick]  (4) -- (12); 

\draw[-,cerulean,ultra thick]  (9) -- (12); 
\draw[-,cerulean,ultra thick]  (9)-- (2); 
\coordinate (P1) at (0,0,0);
\node[Mark2] at (P1) {};
\coordinate (P2) at (0,0,2*3);
\node[Mark2] at (P2) {};
\coordinate (P3) at (0,2*1,2*3);

\draw[dotted,ultra thick,cerulean] (P1) -- (P2);
\draw[-,cerulean,ultra thick] (P2) -- (9);
\draw[dotted,ultra thick,cerulean] (P1) -- (4);
\draw[dotted, ultra thick,cerulean] (P1) -- (12);
\draw[-,ultra thick,cerulean] (P2) -- (2);

\coordinate (X) at (0,0,2*4);
\coordinate (XX) at (0,0,-1);
\coordinate (Z) at (0, 2.5, 0);

\draw[->] (P2) -- (X) node[anchor=north west] {$x$};
\draw[dotted, ultra thick] (XX) -- (P1);

\draw[->] (12) -- (Z) node[anchor=south] {$z$};
\coordinate (T1) at (0,0,2*3*1/2);
\node[Mark3] at (T1) {};
\coordinate (T2) at (0, 2*3 *3/16 ,2* 3 * 1/2);
\node[Mark3] at (T2) {};
\coordinate (T3) at (2*3* 9/28 , 2*3*3/28, 2*3*1/2);
\node[Mark3] at (T3) {};
\draw[-, ultra thick, brickred] (T1) -- (T2);
\draw[-, ultra thick, brickred] (T2) -- (T3);
\draw[-, ultra thick, brickred] (T3) -- (T1);

\node[anchor=north] at (T1) {$(\alpha, 0,0)$};
\node[anchor=south east] at (T2) {$(\alpha, 0, \frac{1+\alpha}{c})$};
\node[anchor=west] at (T3) {$\left(\alpha, \frac{(1+\alpha)d}{ad+c}, \frac{1+\alpha}{ad+c}\right)$};

\end{tikzpicture}
\caption{}
\end{figure}
Hence, the area of the cross section is $A = \frac{1}{2}\frac{1+\alpha}{c}\frac{(1+\alpha)d}{ad+c}$. Therefore, the volume of this part is \[V = \int_0^1 \frac{1}{2}\frac{d(1+\alpha)^2}{c(ad+c)}\mathrm{d}\alpha = \frac{7d}{6c(ad+c)}.\]
It follows that the $F$-signature function $s_X(c,d)$ is \[s_X(c,d) = \frac{7d}{6c(ad+c)} - \frac{d}{6c(ad+c)} = \frac{d}{c(ad+c)}.\]
\end{case}
Combining all cases, we obtain the $F$-signature function as a piecewise function as follows:
\begin{theorem}\label{theorem:F-signature-formula}
\[ s_X(c,d) = \begin{cases}\frac{6d^2-c^2}{6d^2(ad+c)} & \text{ if $c\leq d$} \\
\frac{6d^3-(2d-c)^3}{6cd^2(ad+c)}&  \text{ if $ d \leq c \leq 2d$}\\
\frac{d}{c(ad+c)} & \text{ if $  2d\leq c$}
\end{cases}\]
\end{theorem}

\section{Hilbert--Kunz multiplicity can disagree with \texorpdfstring{$F$}{F}-signature}\label{sec:tensor}
As an application of \cref{theorem:Trivedi} and \cref{theorem:B}, we present an example of a ring witnessing the conclusion of \cref{theorem:A}. First, we recall the necessary base change results for Hilbert--Kunz multiplicity and $F$-signature.

\begin{theorem}[{\cite[Theorem 5.6]{Yao06}}, {\cite[Proposition 3.9(b)]{Kun76}}]\label{theorem:flat-maps-with-regular-fibers}
    Let $(R,\fm)\hookrightarrow (S,\fn)$ be a flat local ring homomorphism between excellent local $\F_p$-algebras. Suppose that the closed fiber $S/\fm S$ is a regular local ring. Then,
\begin{enumerate}[label=(\alph*)]
    \item we have the equality of $F$-signatures $s(R) = s(S)$; \label{theorem:flat-maps-with-regular-fibers-a}
    \item for all $e\in\N$, we have the equality $$\dfrac{\ell_R(R/\fm^{[p^e]})}{p^{e\dim(R)}}=\dfrac{\ell_S(S/\fn^{[p^e]})}{p^{e\dim(S)}}.$$ In particular, taking limits we have $\ehk(R) = \ehk(S)$. \label{theorem:flat-maps-with-regular-fibers-b}
\end{enumerate}
\end{theorem}

As mentioned in \cref{sec:introduction}, the example appearing in \cref{theorem:A} is the tensor product (over $\overline{\F_p}$) of the section rings of the surfaces $\mathcal{H}_a$ with respect to certain line bundles. There are many valid choices for the parameter $a$ and for the line bundles in question, but we only include one for simplicity.
\begin{proof}[Proof of \cref{theorem:A}]
    Consider the line bundles $\cL_1=3D_1 + D_4$ on $\mathcal{H}_2$ and $\cL_2 = 3D_1+2D_4$ on $\mathcal{H}_3$ both defined over $k = \overline{\F_p}$. Define the section rings
\begin{align*}
    R_1:=&\bigoplus\limits_{n\geq 0} H^0(\mathcal{H}_2,\cL_1^n)\\
    R_2:=&\bigoplus\limits_{n\geq 0} H^0(\mathcal{H}_3,\cL_2^n)
\end{align*}
with homogeneous maximal ideals $\fm_i\subseteq R_i$. Finally, let $R:=R_1\otimes_k R_2$. Since they are toric, the rings $R_i$ are both strongly $F$-regular \cite[\S 6]{Smi00} and hence so too is $R$ by \cite[Theorem 5.2]{Has03}. Set $\fp_i:=\fm_i R$ and let $A_i$ be the localization of $R_i$ at $\fm_i$.
Note that the fiber $Y_1$ over $\fm_1$ of the map $ \Spec(R)\to \Spec(A_1)$ is isomorphic to $\Spec(R_2)$. Since $\Spec(R_2)$ is the affine cone over $\mathcal{H}_3$, we may pick a maximal ideal $\fn_1 \in \Spec(R)$ that lies in the fiber over $\fm_1$ and such that $Y_1$ is regular at $\fn_1$ (we can do this by picking any closed point of $\Spec(R_2)$ other than the homogeneous maximal ideal). In other words, we have a flat local ring homomorphism
    \begin{equation*}
        (A_1,\fm_1 A_1)\to (R_{\fn_1},\fn_1 R_{\fn_1})
    \end{equation*}
    whose closed fiber is isomorphic to the local ring of $Y_1$ at $\fn_1$, which is regular by construction. Similarly, we can pick a maximal ideal $\fn_2$ of $R$ in the fiber over $\fm_2$ of $\Spec(R) \to \Spec(A_2)$ such that the closed fiber of $ (A_2,\fm_2 A_2)\to (R_{\fn_2},\fn_2 R_{\fn_2})$ is regular.
    It follows from \cref{theorem:flat-maps-with-regular-fibers} that $s(A_i) = s(R_{\fn_i})$ and $\ehk(A_i) = \ehk(R_{\fn_i})$ for $i=1,2$. Using the formulas of \cref{theorem:Trivedi,theorem:F-signature-formula} we obtain
\[
\begin{array}{rclclcl}
\ehk(R_{\fn_1}) &=& \frac{188}{45} &<& \frac{139}{18} &= & \ehk(R_{\fn_2}), \\
s(R_{\fn_1})    &=& \frac{1}{15}   &<& \frac{47}{648} &= & s(R_{\fn_2})
\end{array}
\]
as desired.
\end{proof}

\printbibliography

\end{document}